\documentclass[a4paper]{amsart}

\usepackage[english]{babel} % (british english
\usepackage[utf8]{inputenc}
\usepackage{amsmath,amssymb}
\usepackage{amsthm}
\usepackage{graphicx}

\usepackage[tableposition=top]{caption}
\usepackage{subcaption}
\usepackage{caption,mathrsfs,longtable}
\usepackage[font={color=black,it},figurename=Figure,labelfont={it}]{caption}

\usepackage{xcolor}
\definecolor{bblue}{rgb}{0.2, 0.4, 0.8}
\definecolor{bgreen}{rgb}{0.2, 0.6, 0.4}
\definecolor{bred}{rgb}{0.8, 0.4, 0.2}
\definecolor{bviolet}{rgb}{0.7, 0.2, 0.7}
\definecolor{blackred}{rgb}{0.6, 0.3, 0.3}
\definecolor{blackblue}{rgb}{0.3, 0.3, 0.6}

\usepackage{hyperref}
\definecolor{lgreen}{rgb}{0.0, 0.48, 0.0}
\definecolor{lpurple}{rgb}{0.48, 0.0, 0.48}
\hypersetup{linktocpage,
            colorlinks=true,
            linkcolor=lgreen,
            citecolor=lpurple,
            linktoc=true}
\usepackage[nameinlink]{cleveref} % I want to reference a range of figures

\usepackage[nocompress]{cite}

\newtheorem{theorem}{Theorem}
\newtheorem{proposition}[theorem]{Proposition}

\newtheorem{definition}[theorem]{Definition}

\renewcommand\paragraph\subsection

\DeclareMathOperator{\Poiss}{Poisson}
\DeclareMathOperator{\Geom}{Geometric}
\DeclareMathOperator{\Var}{\mathbb{V}}
\DeclareMathOperator{\Beta}{Beta}
\DeclareMathOperator{\Uniform}{Uniform}
\newcommand{\convergesto}{\overset{d}{\longrightarrow}}
\def\leq{\leqslant}
\def\geq{\geqslant}

\usepackage{multirow}

\usepackage[foot]{amsaddr}

\newcommand{\mylongtitle}[1]{%
  \ifodd\value{page}%
    \protect\parbox{0.97\linewidth}{#1}\hfill%
  \else%
    \hfill\protect\parbox{0.97\linewidth}{#1}%
  \fi%
}

\title[\mylongtitle{Asymptotic Distribution of Parameters in Random Maps}]
{Asymptotic Distribution of\\ Parameters in Random Maps}

\author{
    Olivier Bodini$^{1}$
}
\address[1]{
    UMR CNRS 7030, LIPN,
    Universit\'{e} Paris 13,
    99 Avenue Jean Baptiste Cl\'{e}ment, 93430 Villetaneuse, FRANCE.
}
\thanks{{This work was partially supported by the METACONC
Project of ANR-MOST (2016--2019).}}
\email{\texttt{\{olivier.bodini, dovgal\}-at-lipn.fr}}
\author{
    Julien Courtiel$^{2}$
}
\address[2]{
    UMR 6072, AMACC,
    Universit\'{e} de Caen Normandie
    Campus C\^{o}te de Nacre, 
    Boulevard du Mar\'{e}chal Juin
    14032 Caen, FRANCE.
}
\email{\texttt{julien.courtiel-at-unicaen.fr}}
\author{
    Sergey Dovgal$^{1,3,4}$
}
\address[3]{
    IRIF,
    Universit\'{e} Paris 7, 5 Rue Thomas Mann
    75013 Paris, FRANCE
}
\address[4]{
    Moscow Institute of Physics and Technology, 
    Institutskiy per. 9, Dolgoprudny,
    141700 RUSSIA
}
\author{
    Hsien-Kuei Hwang$^{5}$
}
\address[5]{
    Institute of Statistical Science, 
    Academia Sinica,
    128, Section 2, Academia Rd, 
    Nangang District, Taipei City, 11529 TAIWAN.
}
\email{\texttt{hkhwang-at-stat.sinica.edu.tw}}

\graphicspath{{./images/}}

\begin{document}

\maketitle

\begin{abstract}
We consider random rooted maps without regard to their genus, with
fixed large number of edges, and address the problem of limiting
distributions for six different parameters: vertices, leaves, loops,
root edges, root isthmus, and root vertex degree. Each of these leads
to a different limiting distribution, varying from (discrete)
geometric and Poisson distributions to different continuous ones:
Beta, normal, uniform, and an unusual distribution whose moments are
characterised by a recursive triangular array.

\smallskip
\noindent \textbf{\keywordsname.}
{
\renewcommand\and{\!\!, }
\keywords{Random maps \and Analytic combinatorics \and Rooted Maps
\and Beta law \and Limit laws \and Patterns \and Generating functions
\and Riccati equation.}
}
\end{abstract}

\section{Introduction}
\label{section:introduction}

\subsection{Motivation for our work}
Rooted maps form a ubiquitous family of combinatorial objects, of
considerable importance in combinatorics, in theoretical physics, and
in image processing. They describe the possible ways to embed graphs
into compact oriented surfaces \cite{LZgraphs}.

The present paper focuses on asymptotic enumeration of basic 
parameters in rooted maps with no restriction on genus. From a
generating function viewpoint, if the genus of the maps is not fixed,
then the generating function of rooted maps is non-analytic (namely,
convergent only at zero) and often satisfies a Riccati differential
equation, in contrast to \emph{planar maps} for which analytic
(convergent) generating functions abound. The divergent Riccati
equations appear frequently in enumerative combinatorics. For
example, at least 39 entries in Sloane's OEIS~\cite{OEIS} were found
containing sequences whose generating functions satisfy Riccati
equations, including some entries related to the families of
indecomposable combinatorial objects, moments of probability
distributions, chord diagrams
\cite{terminalchords,CoYeZe,flajolet2000analytic}, Feynman diagrams
\cite{feynman}, etc. Some of these are closely connected to maps.
Indeed, it is known that rooted maps with no genus restriction also
encode different combinatorial families such as chord diagrams and
Feynman diagrams on the one hand, and different fragments of lambda
calculus \cite{BoGaJa2013,ZGcorr} on the other hand. Thus most
asymptotic information obtained on maps can often be transferred to
the aforementioned objects and lead to a better understanding of them
in the corresponding domains.

While the asymptotics and stochastics on planar maps have been 
extensively studied (see for example 
\cite{banderier2001random,bernardi2017counting,Bender1986,drmota2013central,liskovets1999pattern}), those on 
rooted maps with no genus restriction have received comparatively 
much less attention in the literature. 
%
% Strangely,
% despite a huge literature on the combinatorics of maps, this aspect
% has received little attention in the literature. \todo{\textbf{what
% follows in this paragraph requires a rewrite: There exist certain
% results on the distributions of patterns in planar maps: see e.g. a
% study of the number of vertices of given
% degree~\cite{drmota2013central}. One of the examples of non-Gaussian
% limit law in a combinatorial structure is given
% in~\cite{banderier2001random} with Airy distribution (this paper
% involves planar maps as well). In~\cite{bernardi2017counting},
% differential equations for coloured planar maps are considered.}}
Of closest connection to our study here is the paper by Arquès and
Béraud \cite{arques2000rooted}, which contains several
characterisations of the number of rooted maps and their generating
functions. In particular, they give an explicit formula for the
number of maps, expressed as an infinite sum, from which the
asymptotic number of maps with $n$ edges can be deduced (which is
$(2n+1)!!$). Recently, Carrance~\cite{carrance2017uniform} obtained
the distribution of genus in bipartite random maps. To our knowledge,
no other asymptotic distribution properties of map statistics have 
been properly examined so far. Along a different direction, Flajolet 
and Noy~\cite{flajolet2000analytic} investigated basic statistics on 
chord diagrams, and Courtiel and Yeats \cite{terminalchords} studied
the distribution of \emph{terminal chords}.

From an asymptotic point of view, for planar enumeration, as Bender
and Richmond put it in \cite{Bender1986}: ``\emph{The two most
successful techniques for obtaining asymptotic information from
functional equations of the sort arising in planar enumeration are
Lagrange inversion and the use of contour integration.}'' An equally
useful analytic technique is the saddle-point method as large powers
of generating functions are ubiquitous in map asymptotics; see
\cite{banderier2001random,Flajolet2009} for more detailed
information. In contrast, for divergent series, Odlyzko writes in his
survey \cite{Odlyzko1995}: ``\emph{There are few methods for dealing
with asymptotics of formal power series, at least when compared to
the wealth of techniques available for studying analytic generating
functions.}" We show however that a few simple linearizing techniques
are very helpful in deriving the diverse limit laws mentioned in the
Abstract; the approaches we use may also be of potential application
to other closely related problems.

\subsection{Definitions}
For a rigorous definition of a rooted combinatorial map we refer, for
example, to~\cite{LZgraphs,arques2000rooted}. For our purposes in
this extended abstract we use a less formal but more intuitive
definition.
\begin{definition}[Maps]
    A \emph{map} is a connected multigraph endowed with a cyclic
    ordering of consecutive half-edges incident to each vertex.
    Multiple edges and loops are allowed.
    Around each vertex, each pair of adjacent half-edges is said to 
    form a \emph{corner}.
    If there is only one half-edge, there is only one corner.
    A \emph{rooted map} is a map with a distinguished corner.
\end{definition}

\begin{figure}[hbt!]
    \begin{center}
        \includegraphics[width=0.7\textwidth]{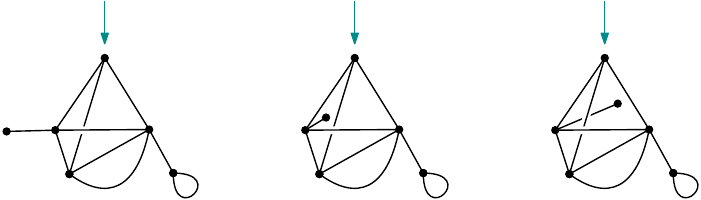}
    \end{center}
    \caption{Three rooted maps. Each root is marked by an arrow.  
    The two last maps are equal.}
    \label{fig:mapexample}
\end{figure}

\begin{figure}[hbt!]
    \begin{center}
        \includegraphics[width=0.7\textwidth]{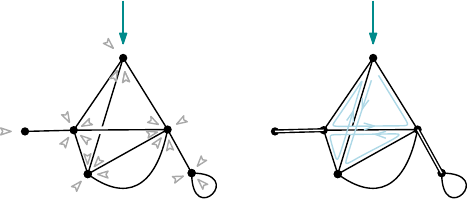}
    \end{center}
    \caption{\textit{Left:} The small triangles point at every corner 
    of the map. \textit{Right:} The light-blue line marks the contour 
    of one face of the map. The double-lined edges are the isthmi of
    the map. The only loop of the map is adjacent to the rightmost
    isthmus, and the vertex incident to this loop has degree $3$.}
    \label{fig:mapexample2}
\end{figure}

\Cref{fig:mapexample} shows some examples of rooted maps. Observe
that the first two maps are different since the cyclic ordering is
not the same: in the first map, the pendant edge follows
counterclockwise the edge after the root (the node pointed to by an
arrow), while in the second map it precedes in counterclockwise
order. In contrast, the last two maps are equal: although the leaves
are at different positions, one can find an isomorphism between the
two maps preserving the vertices, the root and the cyclic orderings
around each vertex. The corners of the leftmost map are displayed
in~\Cref{fig:mapexample2} (left), showing all the possible rootings
of this map.

\begin{definition}[Map features]
    A \emph{face} can be obtained by starting at some corner, moving
    along an incident half-edge, then switching to the next clockwise
    half-edge and repeating the procedure until the starting corner
    is met. A \emph{loop} is an edge that connects the same vertex.
    An \emph{isthmus} is an edge such that the deletion of this edge
    increases the number of connected components of the underlying 
    graph.
    The \emph{degree} of a
    vertex is the number of half-edges incident to this vertex.
\end{definition}

These definitions are illustrated in~\Cref{fig:mapexample2} (right).

Arqu\`{e}s and B\'{e}raud~\cite{arques2000rooted} prove that the
generating function of maps \( M(z) := \sum_{n \geq 0} m_n z^n \),
where \( m_n \) enumerates the number of maps with \( n \) edges,
satisfies\
\begin{equation}
    \label{eq:riccati:maps}
    2z^2 M'(z) = (1-z)M(z) - 1 - z M(z)^2 ,
\end{equation}
a typical Riccati equation whose first few Taylor coefficients read 
$M(z) = 1+2z+20z^2+444z^3 +16944z^4+\cdots$.

\begin{table}[!hbt]
\caption{The six map statistics and their limit laws studied in this 
extended abstract.}
\label{table:comparison}
%-- macro required for proper horizontal alignment
\newcommand\ok[1]{\rlap{$\;#1$}{\hphantom{M}}}
\begin{tabular}{lp{5cm}ll}
    % \multicolumn{1}{l}{\textbf{GF}} &
	\multicolumn{1}{l}{\textbf{Statistics}} &
	\multicolumn{1}{l}{\textbf{Differential equation}} &
    \multicolumn{1}{l}{\textbf{Mean}} &
	\multicolumn{1}{l}{\textbf{Limit law}} \\ \hline \\[-10pt] 
    % \textbf{GF} & \textbf{Pattern} &
    % \centering \textbf{Differential equation}
    % & \textbf{Limit law} \\[2pt]
    % \hline \\[-2mm]
    % $M$ & maps (edges)       &
    % \( M = 1 + zM + zM^2 + 2 z^2 \partial_z M \)
    % & &  \\[3pt]
    % $L$  &
    \multirow{2}{*}{leaves} & 
    \( \ok{L} = v + (2-u)z L + z L^2 \)
    & \multirow{2}{*}{$1$} &
    \multirow{2}{*}{\( \Poiss(1) \)}
    \\
    & \( 
    \qquad { } + 2 z^2 \partial_z L 
    +  z(1-v) \partial_v L \)
    & &
    \\[3pt]
    % $C$ &
    root isthmic parts   &
    \( \ok{C} = 1 + zC + vz  C|_{v = 1}C + 2z^2 \partial_z C \)
    & $2$ &
    \( \Geom(\frac12) \)
    \\[4pt]
    % $X$ &
    vertices  &
    \( \ok{X} = v + z X + zX^2 + 2z^2 \partial_z X \)
    & $\log n$ &
    $\mathcal N(\log n, \log n)$
    \\[3pt]
    % \multirow{2}{*}{$Y$} &
    \multirow{2}{*}{loops}
    &
    \( \ok{Y} = v + vzY + vz Y|_{v = 1} Y   \)
    & \multirow{2}{*}{$\frac12n$} &
    \multirow{2}{*}{\(\text{A new law}^* \)}
    \\
    &
    \(  \qquad { } + 2vz^2  \partial_z Y + v^2z(vw-1) 
    \partial_v Y \)
    &\\ [3pt]
    % $E$ &
    \multirow{2}{*}{root edges}         &
    \( \ok{E} = 1 + vz E + vz E|_{v = 1}E \)
    & \multirow{2}{*}{$\frac23n$} &
    \multirow{2}{*}{\(  \Beta(1, \frac12) \)}
    \\
    & 
    \( \qquad {} + 2vz^2 \partial_z E \)
    & & 
    \\[3pt]
    % \multirow{2}{*}{$D$} &
    \multirow{2}{*}{root degree}
    &
    \( \ok{D} = 1 + v^2z D + vz D|_{v=1} D 
     \)
     & \multirow{2}{*}{$n$} &
    \multirow{2}{*}{\(\text{Uniform}[0,2]\)}
    \\
    & 
    \( \qquad {} + 2vz^2 \partial_z D 
    - v^2(1-v)z \partial_v D \)
    &
    \\%[3pt]

\end{tabular}
\end{table}

\subsection{Results and methods}

We address in this paper the analysis of the extended equations of
\eqref{eq:riccati:maps} for bivariate (and in one case, trivariate)
generating functions \( M(z, v) := \sum_{n, k \geq 0} m_{n,k} z^n v^k
\), where \( m_{n,k} \) stands for the number of maps with \( n \)
edges and the value of the shape parameter equal to \( k \). We
obtain limit laws for the distributions of six different
parameters~(see~\Cref{fig:root-vertex-pgf,fig:root-isthmus-pgf,fig:vertices-pgf,%
fig:root-edges-pgf,fig:joint,fig:loops-pgf}).

\begin{figure}[!htb]
    \begin{minipage}{.45\textwidth}
        \centering
        \includegraphics[width=\textwidth]{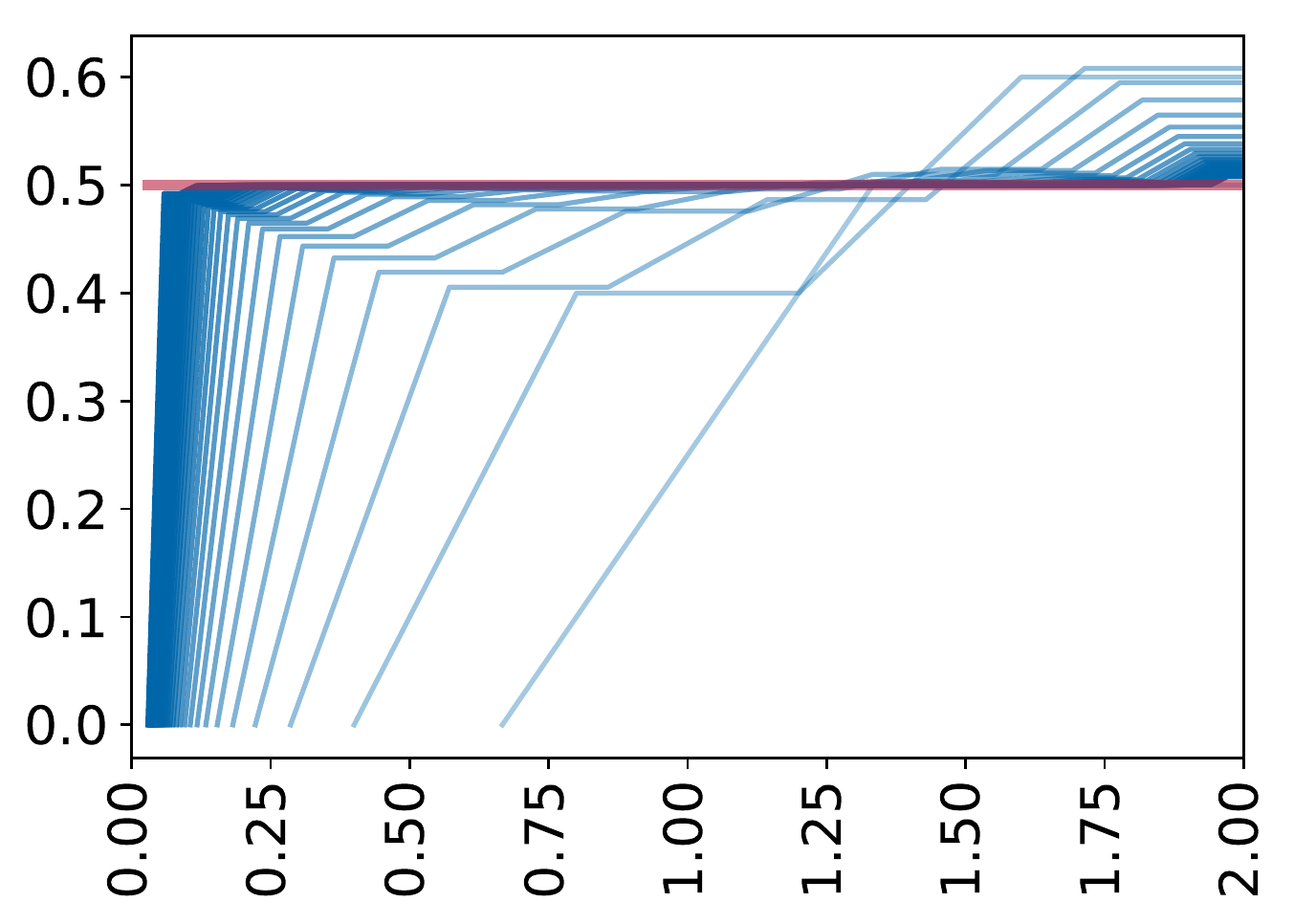}
    \end{minipage}%
    \hfill
    \begin{minipage}{.45\textwidth}
        \centering
        \includegraphics[width=\textwidth]{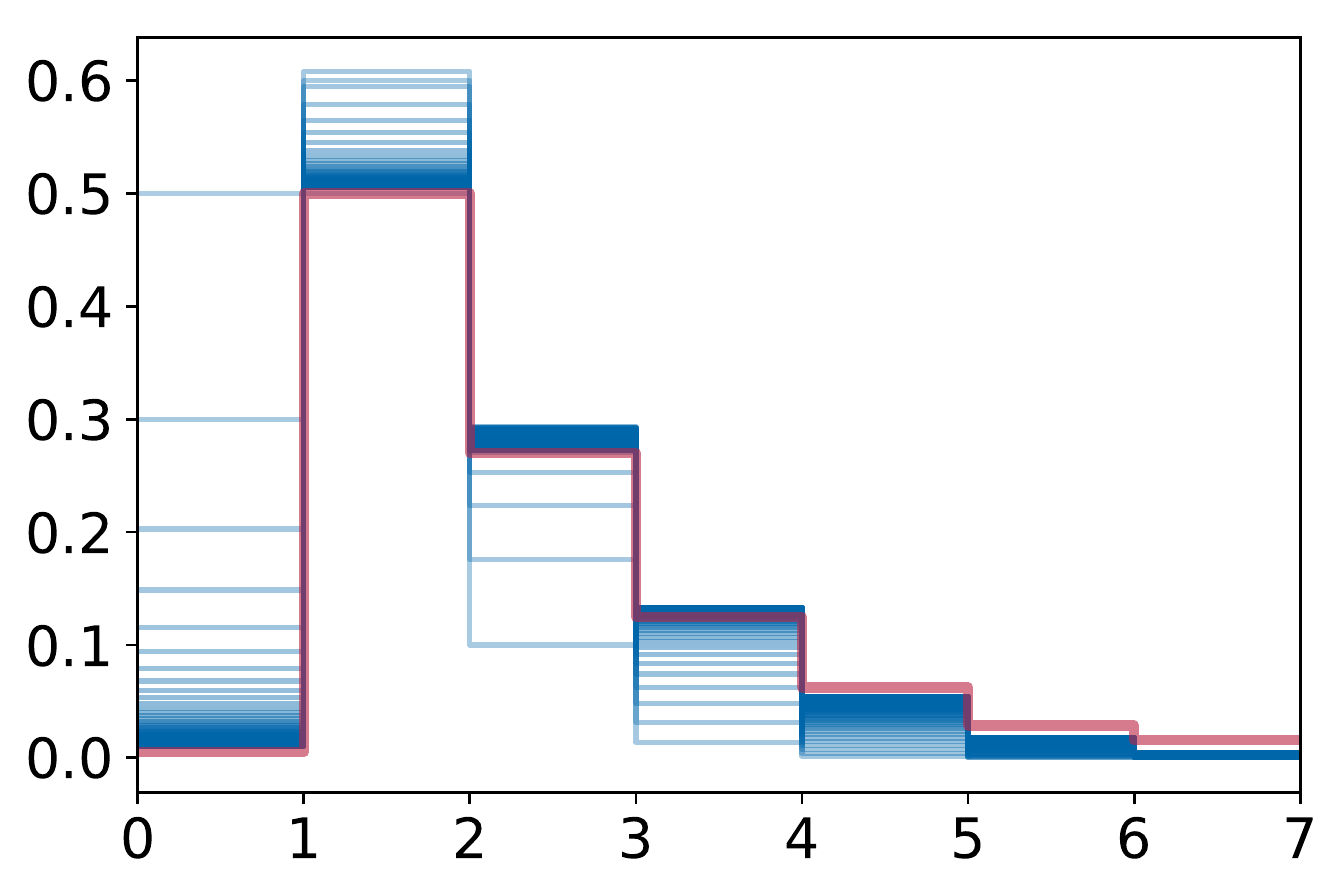}
    \end{minipage}
        \caption{\textit{Left: }Root vertex degree. 
        \textit{Right:} Number of root isthmic parts.}
        \label{fig:root-isthmus-pgf}
        \label{fig:root-vertex-pgf}
\end{figure}                   

\begin{figure}[!htb]
    \begin{minipage}{.45\textwidth}
        \centering
        \includegraphics[width=\textwidth]{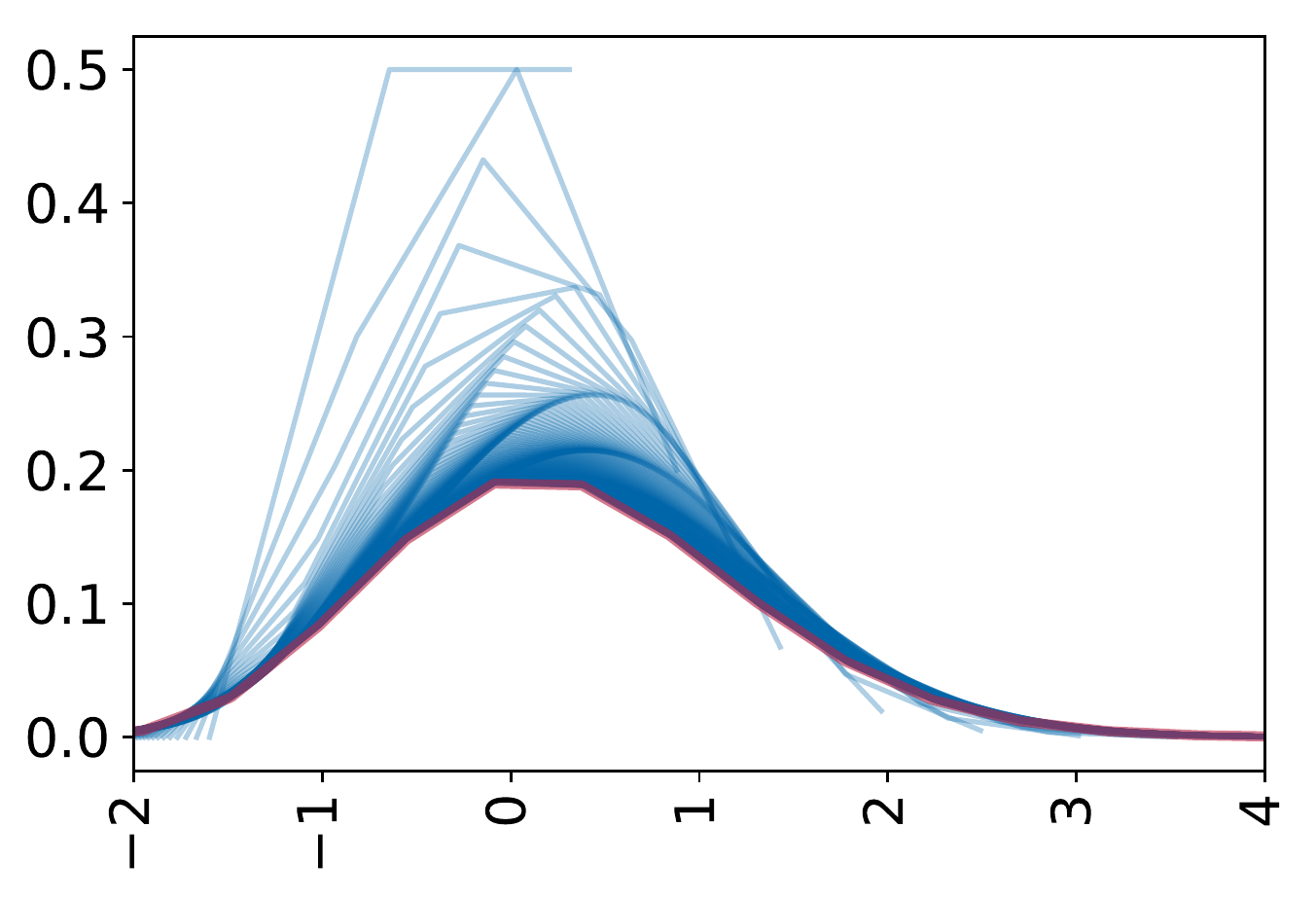}
    \end{minipage}
    \hfill
    \begin{minipage}{.45\textwidth}
        \centering
        \includegraphics[width=\textwidth]{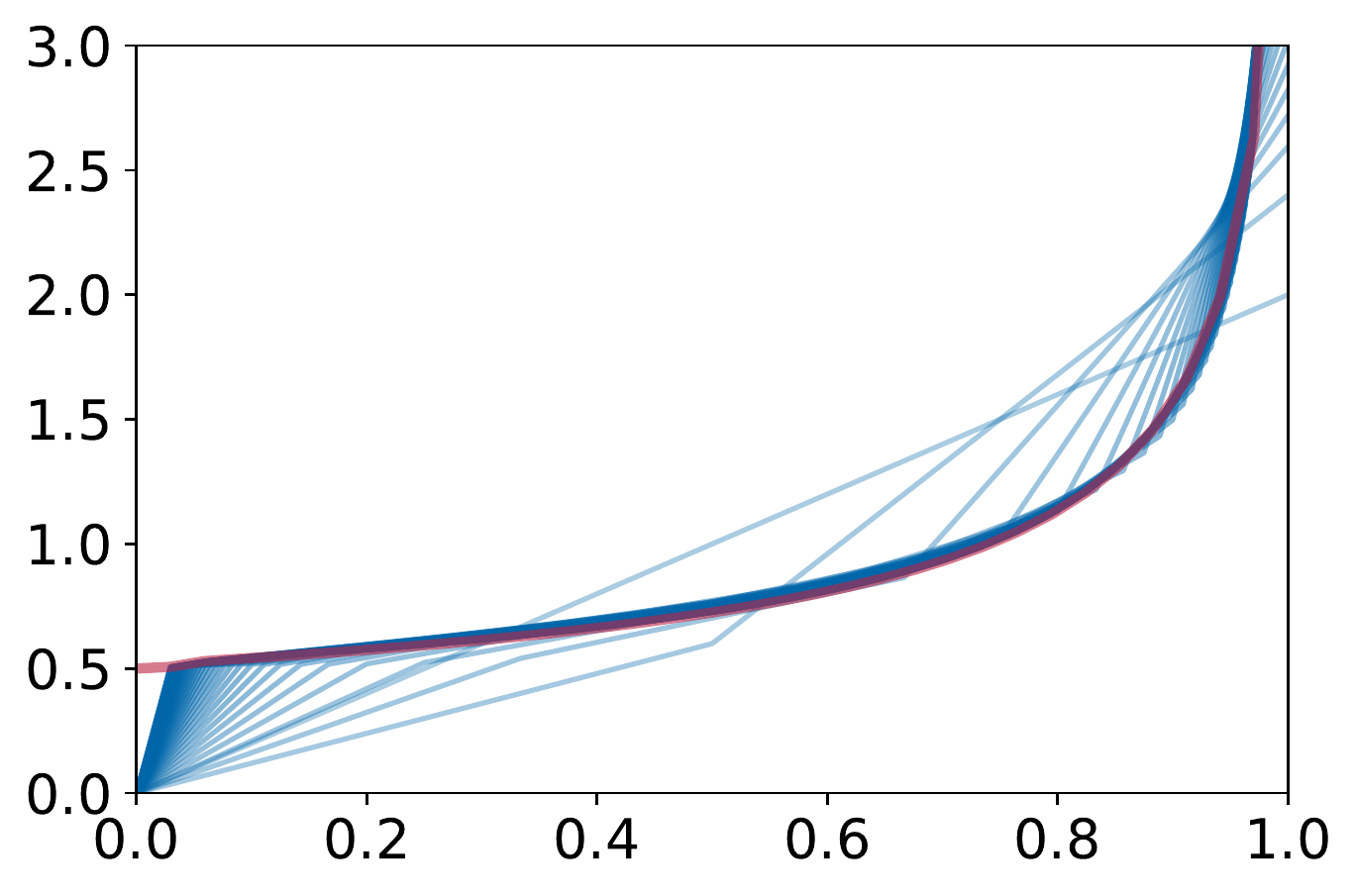}
    \end{minipage}
    \caption{\textit{Left:} Number of vertices. 
    \textit{Right:} Number of root edges.}
        \label{fig:vertices-pgf}
    \label{fig:root-edges-pgf}
\end{figure}

\begin{figure}[!htb]
    \begin{minipage}{.55\textwidth}
        \includegraphics[width=\textwidth]{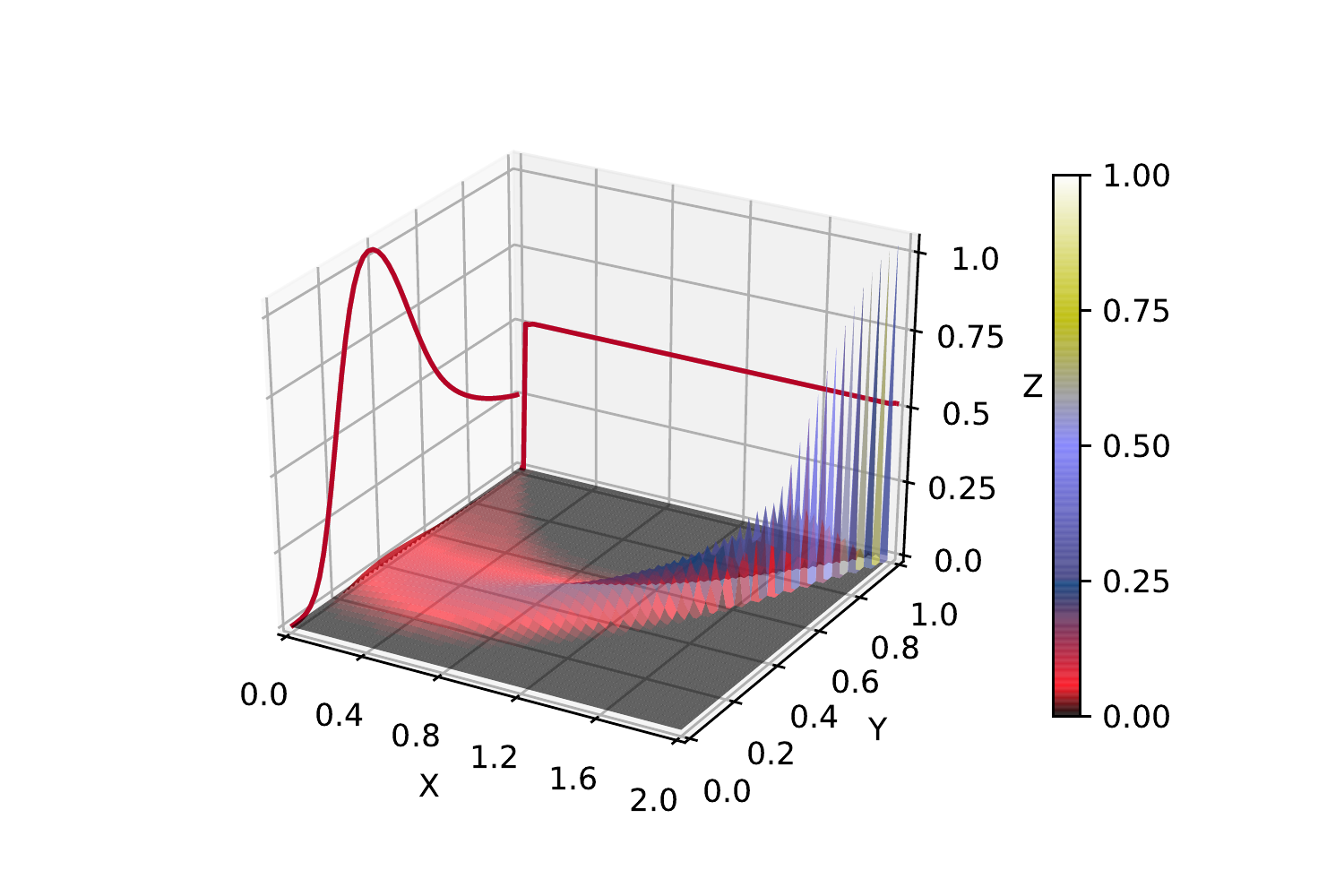}
    \end{minipage}\hfill
    \begin{minipage}{0.45\textwidth}
        \centering
        \includegraphics[width=\textwidth]{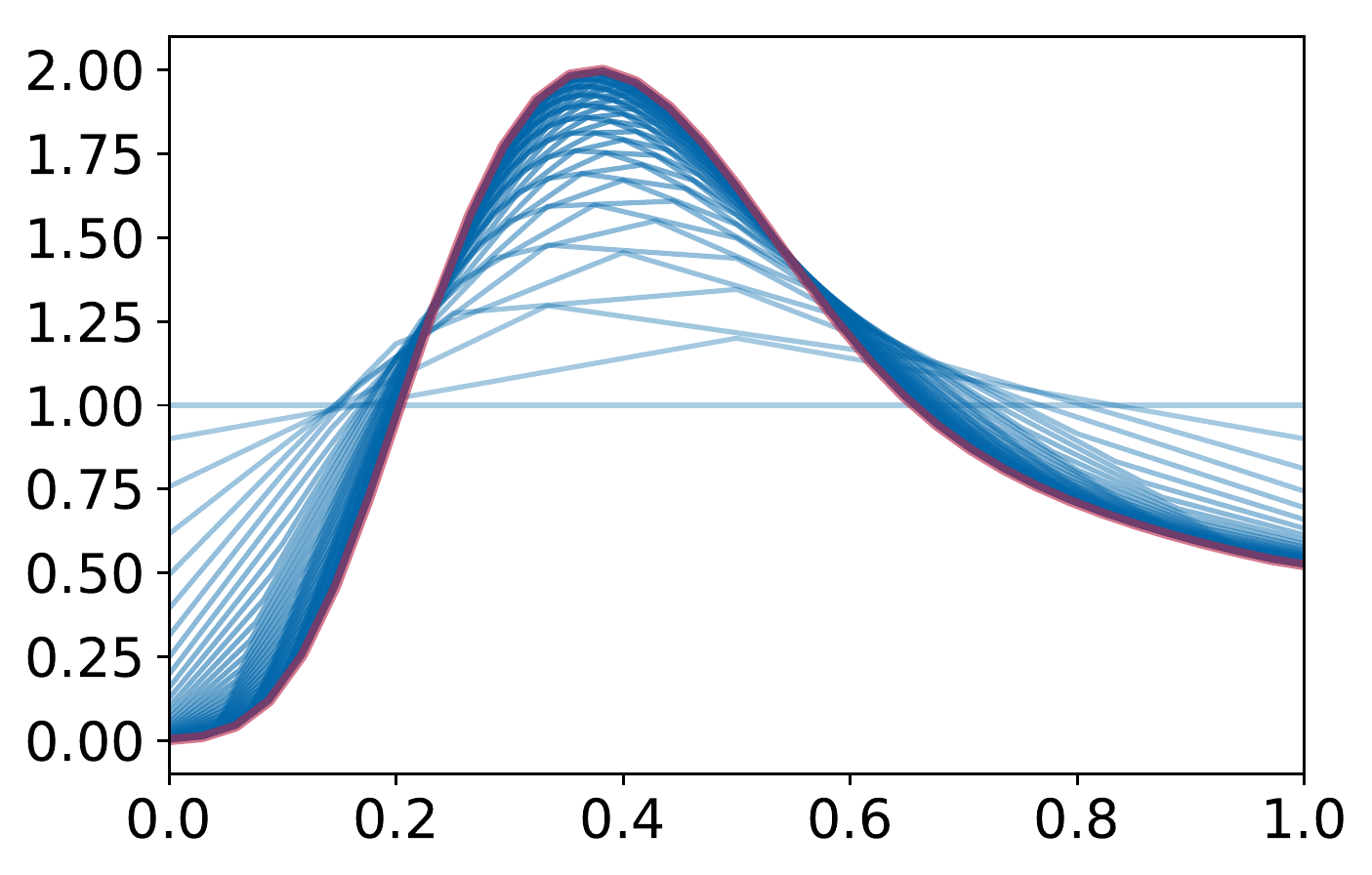}
    \end{minipage}%
     \caption{\textit{Left:} Joint distribution of root vertex 
     degree and the number of loops. \textit{Right:} Number of loops.}
          \label{fig:loops-pgf}  
          \label{fig:joint}              
\end{figure}

We collect the statistics and their limit laws studied here
in~\Cref{table:comparison} for comparison. We see that some of the
limit laws are discrete (Poisson and Geometric), one of them (the
number of vertices) is Gaussian with a logarithmic mean, which we 
denote by $\mathcal{N}(\log n,\log n)$, and the others are
continuous. For the number of root edges, root degree and loops, the
corresponding limit laws are normalized by \( n \), the total number
of edges. The distribution of the number of loops follows a rather
unusual limit law (see \Cref{fig:loops-pgf}) in the sense that we can 
only characterise the limit law by its moment sequence, $\eta^{}_l$,
which satisfies $\eta^{}_l=\eta^{}_{0,l}$ with $\eta^{}_{k,l}$
computable only through a recurrence involving $\eta^{}_{k-1,l}$ and
$\eta^{}_{k+1,l-1}$. The corresponding probability density function
of this law remains unknown and does not have an explicit expression
at this stage (see~\Cref{fig:joint}). Finally, by the bijection
from~\cite{CoYeZe} and a known property of chord diagrams
in~\cite{flajolet2000analytic}, it is possible to deduce the limit
laws for the number of leaves.

One technique we use several times in our proofs consists in
linearising the differential equations satisfied by the generating
functions, by choosing a suitable transformation, inspired from the
resolution of Riccati equations. Once the dominant term is
identified, the analysis for the limit law becomes more or less
straightforward. When such a technique fails, we rely then on the
method of moments, which establishes weak convergence by computing
all higher derivatives of $M(z,v)$ at $v=1$ and by examining
asymptotically the ratios $[z^n]\partial_v^k
M(z,v)|_{v=1}/[z^n]M(z,1)$ (which correspond to the moments of random
variable). Such a procedure also linearises to some extent the more
complicated bivariate nature of the differential equations and
facilitates the resolution complexity of the asymptotic problem.

\paragraph{Structure of the Paper.}
In~\Cref{section:differential:equations} we derive the nonlinear 
differential equations satisfied by the generating functions of the 
map statistics. Then in~\Cref{section:limit:laws} we sketch the 
proofs for the limit laws of five statistics based on generating 
functions. The Poisson law for the number of leaves (together with 
the root face degree and the number of trivial loops) will be proved 
by a direct combinatorial approach in the last section. 

%%%%%%%%%%%%%%%%%%%%%%%%%%%%%%%%%%%%%%%%%%%%%%%%%%%%%%%%
\section{Differential equations for maps}
\label{section:differential:equations} 

In this section, we derive the differential equations satisfied by 
the bivariate or trivariate generating functions with the additional 
variable(s) counting the shape statistics.

\begin{figure}[hbt!]
    \begin{center}
    \includegraphics[width=0.6\textwidth]{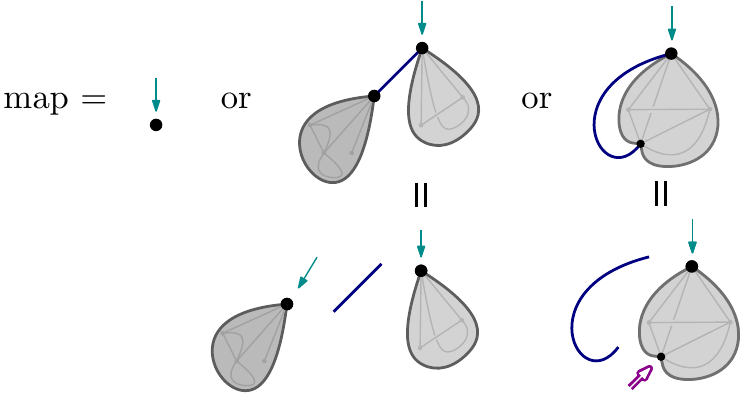}
    \end{center}
    \caption{A symbolic construction of rooted maps.}
    \label{fig:symbolic:maps}
\end{figure}

\paragraph{Univariate generating function of maps.}
Since the Riccati equation \eqref{eq:riccati:maps} lies at the basis
of all other extended equations in~\Cref{table:comparison}, we give a
quick proof of it via the recurrence satisfied by $m_n$, the number
of maps with \( n \) edges (see~\Cref{fig:symbolic:maps}):
\begin{equation}
    \label{eq:recurrence:maps}
    m_n 
    = \boldsymbol 1_{[n=0]} 
    + \sum_{0\leq k<n} m_k m_{n-1-k} 
    + (2n-1) m_{n-1},
\end{equation}
which then implies the Riccati equation \eqref{eq:riccati:maps}.

First, \( m_0 = 1 \) because there is only one map with \( 0 \)
edges. Then a map with \( n \) edges can be formed either by
connecting the roots of two maps (with \( k \) and \( n-k-1 \) edges,
respectively) with an isthmus, or by adding an edge to a map with \(
n - 1 \) edges, connecting the root and a corner. The number of
possible ways to insert an edge in this way is equal to \( 2 n-1 \),
because there are $2n-2$ corners in a map of size $n-1$, and there
are two possible ways to insert a new edge at the root corner (either
before, or after the root). This proves \eqref{eq:recurrence:maps}.

\paragraph{Vertices.}
Consider now the bivariate generating function
\(
    X(z, v) = \sum_{n, k \geq 0} x_{n,k} z^n v^k,
\)
where \( x_{n,k} \) is equal to the number of rooted maps with \( n
\) edges and \( k \) vertices. Arqu\`{e}s and
B\'{e}raud~\cite{arques2000rooted} showed that
\begin{equation}
    \label{eq:de:vertices}
    X(z, v) = v
    + z X(z, v)    
    + z X(z,v)^2
    + 2 z^2 \partial_z X(z, v).
\end{equation}
This recurrence can be obtained from~\eqref{eq:recurrence:maps} by
noticing that no new vertex is created when we connect two maps with
an isthmus, nor when we add a new root edge to a map. Note that \(
X(z, v) \) satisfies another functional equation
(see~\cite{arques2000rooted})
\[
    X(z, v) = v + z X(z, v) X(z, v+1),
\]
which seems less useful from an asymptotic point of view. 

\paragraph{Root isthmic parts.}
We count here the \emph{root isthmic parts}, which are the number of
isthmic constructions used at the root vertex. Note that an isthmic
part may not be a bridge because the additional edge constructor may
induce additional connections.

We show that the bivariate generating function \( C(z, v) = \sum_{n,
k \geq 0} c_{n,k} z^n v^k \), where \( c_{n,k} \) enumerates the 
number of maps with \( n \) edges and \( k \) root isthmic parts, 
satisfies 
\begin{equation} \label{eq:Czu}
    C(z, v) = 1 + z C(z, v) + v z C(z, v) C(z, 1)
    + 2 z^2 \partial_z C (z, v).
\end{equation}

In~\Cref{fig:symbolic:maps}, the number of root isthmic parts only
changes whenever two maps are connected by an isthmus. This yields \(
vz C(z, v) C(z, 1) \) instead of \( z C^2 \).

\paragraph{Root edges.}
Similarly, consider \( E(z, v) = \sum_{n, k \geq 0} e_{n,k} z^n v^k
\), where \( e_{n,k} \) counts the number of rooted maps with \( n \)
edges and \( k \) root edges. We show that \( E(z, v) \) satisfies
\begin{equation}\label{eq:Ezv}
    E = 1 + vz E + v z E|_{v=1} E
    + 2 v z^2 \partial_z E .
\end{equation}
This again results from the recurrence~\eqref{eq:recurrence:maps} and
from~\Cref{fig:symbolic:maps}: the non-root edges come from the
bottom map in the isthmic construction, yielding the term \( vz E(z,
v) E(z, 1) \).

\paragraph{Root Degree.}
Consider the degree of the root vertex. Note that this may be
different from the number of root edges because for the root degree,
each loop edge is counted twice, therefore the degree of the root
vertex varies from \( 0 \) to \( 2n \). By duality, the distribution
of the root face degree is the same as the distribution of the root
vertex degree.

Let \( D(z, v) = \sum_{n,k\geq 0} d_{n,k} z^n v^k \) denote the 
bivariate generating function for maps with variable \( v \) marking 
root degree. Then 
\begin{equation}\label{eq:Dzu}
    D = 1 
    + v^2 z D
    + vz D|_{v=1} D
    + 2 v z^2 \partial_z D
    - v^2(1-v) z \partial_v D.
\end{equation}

\begin{figure}[hbt!]
    \begin{center}
    \includegraphics[width=0.8\textwidth]{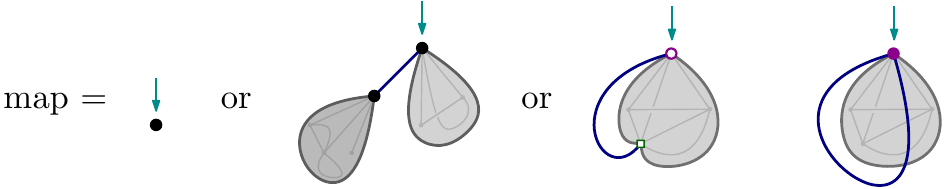}
    \end{center}
    \caption{Symbolic method to count root degree and loops 
    in rooted maps.}
    \label{fig:symbolic:trivial:loops}
\end{figure}

In this case, the original construction in~\Cref{fig:symbolic:maps}
is insufficient, and we need to consider further cases
in~\Cref{fig:symbolic:trivial:loops}. When an additional edge becomes
a loop, it increases the degree of the root vertex by \( 2 \);
otherwise, the root degree is increased merely by $1$. Note that
the equation \eqref{eq:Dzu} is now a \emph{bona fide} partial 
differential equation, making the analysis more difficult. 

\paragraph{Leaves.} The differential equation for the bivariate 
generating function of maps with variable \( v \) marking leaves 
(see~\Cref{table:comparison}) can be obtained in a similar way by
considering different cases in the new edge constructor. The number 
of special leaf corners is equal to the number of leaves.

\paragraph{Loops.}
Finally, we look at the number of loops whose enumeration
necessitates the consideration of the joint distribution of the
number of loops and the number of root edges, namely, we consider 
the trivariate generating function
\(
    Y(z, v, w) = \sum_{n, k, m} y_{n,k,m} z^n v^k w^m
\),
where \( y_{n,k,m} \) denotes the number of rooted maps with \( n \)
edges, root degree equal to \( k \), and \( m \) loops. We show that
\( Y(z, v, w) \) satisfies a partial differential equation
\begin{equation}
    \label{eq:de:loops}
    Y = 1
    + z v Y
    + zv Y|_{v=1} Y
    + 2 z^2 v \partial_z Y
    + z v^2 (vw - 1) \partial_v Y
    .
\end{equation}
As in the symbolic construction of~\Cref{fig:symbolic:trivial:loops},
a new edge becomes a loop only if it is attached to one of the
corners incident to the root vertex. The differential
equation~\eqref{eq:de:loops} is then a modification of~\eqref{eq:Dzu}
with an additional variable counting the number of loops.

Note that Equation~\eqref{eq:de:loops} is \emph{catalytic} with
respect to the variable \( v \), i.e. putting \( v = 1 \) introduces
a new unknown object \( \partial_v Y |_{v=1} \) to the
differential equation. One of the strategies for dealing with
catalytic equations was developed by Bousquet-M\'{e}lou and
Jehanne~\cite{bousquet2006polynomial}, generalising the so-called
kernel method and quadratic method. However, their method does not
work in our case because our equation is differentially algebraic.

%%%%%%%%%%%%%%%%%%%%%%%%%%%%%%%%%%%%%%%%
\section{Limit laws}
\label{section:limit:laws}

This section describes the techniques we employ to establish the 
limit laws.

From now on, by a random map (with $n$ edges) we assume that all
rooted map with $n$ edges are equally likely. For notational
convention, we use $X'=\partial_z X$ to denote derivative with
respect to $z$. Due to space limit, we give only the sketches of the
proofs.

\subsection{Transformation into a linear differential equation}

For most of the equations in the previous section, it turns out that
a transformation similar to that used for Riccati equations largely
simplifies the resolution and leads to solvable recurrences, which
are then suitable for our asymptotic purposes. We begin by solving
the standard Riccati equation \eqref{eq:riccati:maps} and see how a
similar idea extends to other differential equations.

\begin{proposition}
    \label{prop:asymptotic:maps}
    The number \( m_n \) of maps with \( n \) edges satisfies
    \begin{equation}\label{eq:mn-phin}
        \frac{m_n}{\phi_n} = 2n -1 + O\bigl(n^{-1}\bigr), 
        \quad\textrm{where} \quad
        \phi_n = \dfrac{(2n)!}{2^n n!}
        = (2n-1) !!.
    \end{equation}
\end{proposition}
\begin{proof} 
We solve the Riccati equation \eqref{eq:riccati:maps} by considering
the transformation
\begin{equation}
    \label{eq:substitution:maps}
    M(z) = 1 + \dfrac{2 z \phi'(z)}{\phi(z)},
\end{equation}
for some function $\phi(z)$ with $\phi(0)=1$. Substituting this form 
into the equation~\eqref{eq:recurrence:maps}, we get the second-order
differential equation $2z^2 \phi'' + (5z-1)\phi' + \phi =0$. From 
this equation, the coefficients \( \phi_n := [z^n]\phi(z)\) satisfy 
the recurrence
\(
    \phi_{n+1}     
    =
    (2n+1) \phi_n
\), which implies the double factorial form of \( \phi_n \) by 
$\phi_0=1$. 

Moreover, by extracting the coefficient of $z^n$ in
\eqref{eq:substitution:maps}, we obtain a relation between the
coefficients $m_k$ and $\phi_\ell$. By the inequality $m_n \geq
(2n-1) m_{n-1}$ (see \eqref{eq:recurrence:maps}), we then deduce the
asymptotic relation \eqref{eq:mn-phin}. 
\end{proof}

\begin{theorem}
\label{thm:vertices}
    Let \( X_n \) denote the number of vertices in a random 
    rooted map with \( n \) edges. Then $X_n$ follows a central limit 
    theorem with logarithmic mean and logarithmic variance:
\end{theorem}
    \begin{equation}\label{eq:Xn-CLT}
        \dfrac{X_n - \mathbb{E}(X_n)}{\sqrt{\Var(X_n)}} 
        \convergesto \mathcal N(0, 1)
        , \quad
        \mathbb{E}(X_n) \sim \log n
        , \quad
        \Var(X_n) \sim \log n
        .
    \end{equation}
\begin{proof}
Similar to \eqref{eq:substitution:maps}, we define a bivariate 
generating function \( S(z, v) = \sum_{n \geq 0} s_n(v) z^n \) such 
that
\begin{equation*}
    X(z, v) = v + \dfrac{2 z S'}{S}, \quad
    S(0) = 1.                                    
\end{equation*}
Substituting this \( X(z, v) \) into~\eqref{eq:de:vertices} 
%\begin{equation}
%    4z^3 S'' - 2z(1-(3-2v)z) S' + v(1+v)z S = 0
%    .
%\end{equation}
leads to a linear differential equation from which one can extract the
recurrence 
\begin{equation*}
    s_n(v) = \dfrac{(2n+v-2)(2n+v-1)}{2n}\,  s_{n-1}(v).
\end{equation*}
We then get an explicit expression for $s_n(v)$, from which we 
deduce, by singularity analysis, that
\[
    \mathbb{E}(v^{X_n})  
    = \dfrac{2^{v-1}}{\Gamma(v)} \,n^{v-1}
    (1 + O(n^{-1})), 
\] 
and conclude  by applying the Quasi-Powers 
Theorem~\cite{Flajolet2009,hwang1998convergence}. 
\end{proof}

A finer Poisson$(\log n+c)$ approximation, for a suitably chosen $c$, 
is also possible, which results in a better convergence rate $O(\log 
n)^{-1}$ instead of $(\log n)^{-\frac12}$; see \cite{Hwang1999} for 
details. 

\begin{theorem}
    Let \( C_n \) denote the number of root isthmic parts
    in a random rooted map with \( n \) edges. Then, 
    \begin{equation*}
        C_n
        \convergesto
        \Geom\bigl(
            \tfrac12
        \bigr).
    \end{equation*}
\end{theorem}
\begin{proof}
Since \( C(1,z)=M(z) \), we use again the
substitution~\eqref{eq:substitution:maps} and apply it to
\eqref{eq:Czu}:
\begin{equation*}   
    2 z^2 (\phi C' + v \phi' C)
    = (1-(1 + v) z )\phi C - \phi.
\end{equation*}
The trick here is to multiply both sides by \( \phi(z)^{v-1} \)  
and set \( Q(z, v) = \phi(z)^v C(z, v) \). %Then $Q(0)=0$ and 
We then obtain
\begin{equation*}
    2z^2Q' = (1-(1+v)z)Q - \phi^v.
\end{equation*}
Using the recurrence for the normalised coefficients \( \widehat
q_n(v) := q_n(v) / \phi_n \) and dominant-term approximations, we
find that the $n$-th coefficient of $Q$ is proportional to
\[
    \widehat q_n(v) = 
    \dfrac{v}{2n} \sum_{1\leq k\leq n} \left(
        \dfrac{n}{k}
    \right)^{v/2}
    + O(n^{-1/2})
    = 
    \dfrac{v}{2 - v} + O( n^{-\frac12} ).
\]
This corresponds to a (shifted by $1$) geometric distribution with 
parameter \( \frac12 \). By the definition \( Q(z, v) = \phi(z)^v 
C(z, v) \), we deduce that the limiting distribution of \( C_n \) is 
also geometric with parameter \( \frac12 \). 
\end{proof}

\begin{theorem}
    Let \( E_n \) denote the number of edges incident to the root
    vertex in a random rooted map with \( n \) edges. Then $E_n$ 
    follows asymptotically a Beta distribution:
    \begin{equation}\label{eq:beta}
        \dfrac{E_n}{n}
        \convergesto
        \Beta\bigl(1, \tfrac12\bigr),
    \end{equation}
    with the density function 
    \( \frac12 (1 - t)^{-\frac12} \) for \( t \in [0, 1) \).
\end{theorem}
\begin{proof} We use again the substitution \( E(z, 1) = M(z) = 1 + 
2z \frac{\phi'}{\phi}\) in \eqref{eq:Ezv}, giving
\[
    2vz^2(\phi E' + \phi'E)
    = (1-2vz)\phi E - \phi.
\]    
With \( Q(z, v) = \phi(z) E(z, v) \), we then obtain
\begin{equation}
    2vz^2 Q' = (1 - 2 vz) Q - \phi.
\end{equation}
This linear differential equation translates into a recurrence 
for the coefficients \( q_n(v) \) of \( Q(z, v) \), which yields the 
closed-form expression
\begin{equation}
    q_n(v) = 2^n n! \sum_{0\leq j\leq n}
    {2j \choose j} 4^{-j} v^{n-j}.
\end{equation}
Returning to  \( E(z, v) \), we see that its coefficients 
behave asymptotically like \( q_n(v) \). This
implies the Beta limit law \eqref{eq:beta} for the random variable \( 
E_n / n \) since $\binom{2j}{j}4^{-j}\sim (\pi j)^{-1/2}$ for large 
$j$. 
\end{proof}

\begin{theorem} \label{thm:rootdeg}
Let \( D_n \) denote the degree of the root vertex in a random
rooted map with \( n \) edges. Then, \( D_n \), divided by the
number of edges, converges in law to the uniform distribution on
\( [0, 2] \):
\begin{equation}\label{eq:Dn}
    \dfrac{D_n}{n}
    \convergesto
    \Uniform \left[0, 2\right].
\end{equation}
\end{theorem}
\begin{proof}
The substitutions
\begin{equation*}
    D(z, 1) = M(z) = 1 + \dfrac{2z \phi'}{\phi}
   , \quad
    \text{and}
    \quad
    D(z, v) = \dfrac{Q(z, v)}{\phi(z)}
\end{equation*}
lead to a partial differential equation, which in turn yields the
recurrence for the coefficients \( q_n(v) := [z^n]Q(z, v) \):
\[
    q_n(v) = v (2n-1 + v) q_{n-1} 
    - v^2 (1 - v) q_{n-1}'(v) + \phi_n.
\]
We then get the exact solution $q_n(v) = \phi_n (1 + v + \cdots +
v^{2n})$. Accordingly, \( d_n(v) := [z^n] D(z, v)\sim q_n(v) \). This
implies the uniform limit law \eqref{eq:Dn}. 
\end{proof}
A more intuitive interpretation of this uniform limit law is given in 
the next section. 

\subsection{Approximation and method of moments}

Unlike all previous proofs, we use the method of moments to establish
the limiting distribution of the number of loops. The situation is
complicated by the presence of the term involving $\partial_v Y$ in
\eqref{eq:de:loops}, which introduces higher order derivatives with
respect to $v$ at \( v = 1 \) when computing the asymptotic of the 
moments. 
\begin{theorem}
Let \( Y_n \) denote the total number of loops in a random rooted map
with \( n \) edges. Then
\begin{equation}
    \dfrac{Y_n}{n}
    \convergesto
    \mathcal L
    ,
\end{equation}
where \( \mathcal L \) is a continuous law with a computable 
density on \( [0, 1] \).
\end{theorem}
\begin{proof} 
First, we show by induction that there exist constants 
\( \eta^{}_{k, \ell} \), such that as \( n \to \infty \),
\begin{equation}\label{eq:eta-as}
    \left.
    [z^n]
    \partial_v^k \partial_w^\ell
    Y(z, v, w) \right|_{v = w = 1} 
    \sim \eta^{}_{k, \ell} n^{k + \ell + 1} \phi_n
    , \quad
    k, \ell \geq 0
     .
\end{equation}
For \( k = \ell = 0 \) the statement clearly holds.
%   \eqref{eq:mn-asymp} with \( \eta^{}_{0,0} = 2 \).
Let \( y_n^{(k,\ell)}:= 
    \left.
    [z^n]
    \partial_v^k \partial_w^\ell
    Y(z, v, w) \right|_{v = w = 1} 
 \)
for larger \( k, \ell \geq 0 \). By translating~\eqref{eq:de:loops}
into the corresponding recurrence for the coefficients and by
collecting the dominant terms (using the induction hypothesis
\eqref{eq:eta-as}), we deduce that
\begin{equation*}
    y_{n}^{(k, \ell)}
    \sim
      (2n + k) y_{n-1}^{(k, \ell)}
      + \ell y_{n-1}^{(k+1, \ell-1)}
      + (2kn - 2k) y_{n-1}^{(k-1, \ell)}
      + \boldsymbol 1_{[k=0]} y_{n-1}^{(k, \ell)}
       .
\end{equation*}
Accordingly, we are led to the recurrence
\begin{equation*}
    \eta^{}_{k, \ell} = 
    \dfrac{1}{k + 2 \ell + \boldsymbol 1_{[k > 0]} }
    ( 2k \eta^{}_{k-1, \ell} + \ell \eta^{}_{k+1, \ell-1}),
\end{equation*}
for $k+\ell>0$ (provided that we interpret $\eta^{}_{k,\ell}=0$ when
any index becomes negative). In particular, when \( \ell = 0 \), we
obtain the moments of the random variable \( E_n \), the number of
root edges: $\eta^{}_{k, 0} = \frac{2^{k+1}}{k + 1}$, which coincides
with the moments of the uniform random variable \( \Uniform [0,2] \).
Finally, it is not complicated to check that the numbers \(
\eta^{}_{0, \ell} \) satisfy the condition of Hausdorff moment
problem, i.e. \( \eta^{}_{0, \ell} \) uniquely determine the limiting
random variable defined on \( [0, 1] \). 
\end{proof}

\section{Combinatorics of map statistics}
\label{section:combinatorial:approaches}

We examine briefly the combinatorial aspect of the map statistics, 
relying our arguments on the close connection between maps and chord 
diagrams (see \cite{cori2009indecomposable}). 

Recall that a \emph{chord diagram} \cite{flajolet2000analytic} with
\( n \) chords is a set of vertices labelled with the numbers \( \{
1, 2, \ldots, 2n\} \) equipped with a perfect matching. A chord
diagram is \emph{indecomposable} if it cannot be expressed as a
concatenation of two smaller diagrams.

\paragraph{Why the root degree follows a uniform law?} 
We begin with Cori's bijection \cite{cori2009indecomposable} between
rooted maps and indecomposable diagrams. In this bijection, each
chord connecting labels \( i \) and \( j \) corresponds to matching
of the half-edges with labels \( i \) and \( j \). The set of
half-edges incident to each vertex of the resulting map corresponds
to the set of nodes to the right of the starting points of the
so-called \emph{outer chords}, i.e. chords that do not lie under any
other chord.

\begin{proposition} There exists a bijection between rooted maps of
root degree $d$ with $n$ edges, and indecomposable diagrams with
$n+1$ chords such that the vertex $k-2$ is matched with vertex $1$. 
\end{proposition}

Once this proposition is available, it leads to a simpler and more
intuitive proof of \Cref{thm:rootdeg} as follows. In a (not
necessarily indecomposable) diagram, the label of the vertex matched
with $1$ follows exactly a uniform law on $\{2,\dots,2n\}$. But a
diagram is almost surely an indecomposable diagram (because its
cardinality is asymptotically the same); thus the label of the vertex
matched with $1$ divided by $2n$ obeys asymptotically a uniform law
on $[0,1]$ (or Uniform$[0,2]$ if divided by $n$ as in
\Cref{thm:rootdeg}).

\begin{figure}
    \begin{minipage}{.45\textwidth}
        \centering
        \includegraphics[width=\textwidth]{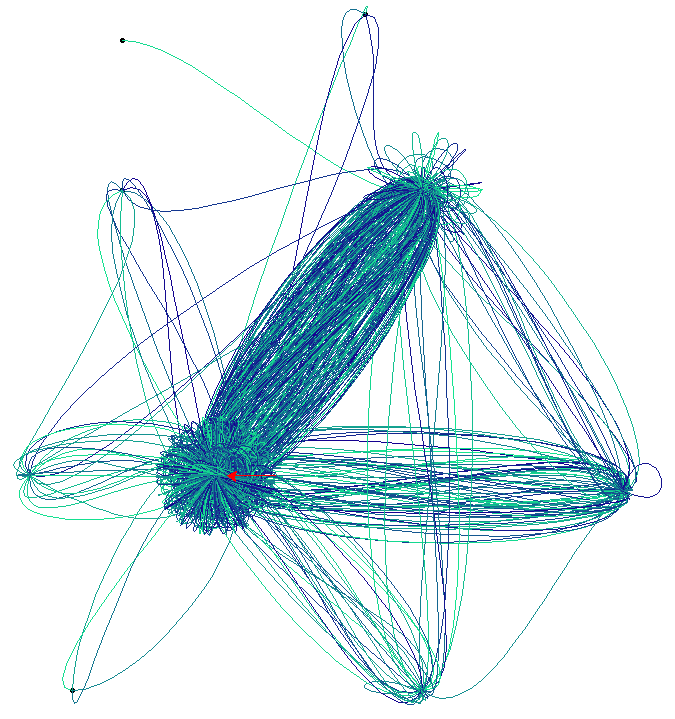}
    \end{minipage}
    \hfill
    \begin{minipage}{.45\textwidth}
        \centering
        \includegraphics[width=\textwidth]{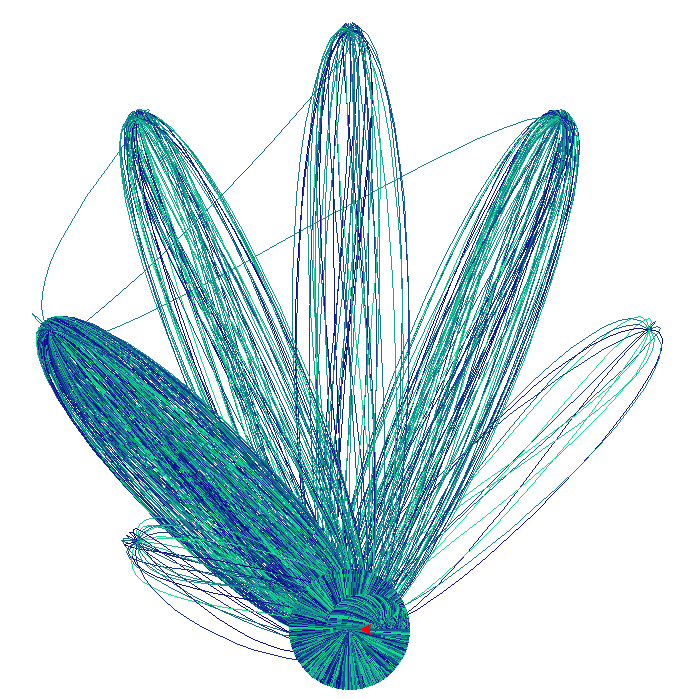}
    \end{minipage}
    \caption{ Random rooted maps, respectively with $1000$ 
    and $20000$ edges.}
    \label{fig:random}
\end{figure}

\paragraph{Uniform random generation.} Cori's bijection is also
useful for generating random rooted maps. Uniformly sampling a random
diagram can be achieved by adding the chords sequentially one after
another. If this procedure results in an indecomposable diagram, it
is rejected (which occurs with asymptotic probability $0$). A  
successful sampled diagram is then transformed into a map using 
Cori's bijection \cite{cori2009indecomposable}. 
Figure~\ref{fig:random} shows two instances of random maps thus 
generated.

\paragraph{The number of leaves.} 
Another bijection in \cite{CoYeZe} is useful in proving the Poisson 
limit law of the number of leaves. This bijection sends leaves of a
map into the \textit{isolated} chords (namely, edges connecting
vertices \( k \) and \( k+1 \)) of an indecomposable chord diagram.
According to~\cite[Theorem 2]{flajolet2000analytic}, the number of
isolated edges in a random chord diagram has a Poisson distribution
with parameter \( 1 \). We can then deduce the following theorem.

\begin{theorem}
The number of leaves in a random map with \( n \) edges follows 
asymptotically a Poisson law with parameter \( 1 \).
\end{theorem}

\paragraph{Two dual parameters.}
We briefly remark that two other parameters, namely \emph{root face
degree} and the number of \emph{trivial loops} do not seem easily
dealt with by the method of generating functions because marking them
requires additional nested information such as the degrees of all the
faces. However, such parameters can be easily marked in their
corresponding dual maps. Their limit distributions are uniform and
Poisson, respectively.

\footnotesize
\bibliographystyle{alpha}
\bibliography{maps-distributions-arxiv}

\end{document}